\documentclass[10pt, leqno]{amsart}
\usepackage{amssymb,amsfonts,amscd,amsbsy}
\usepackage[mathscr]{eucal}
\usepackage{url}

\theoremstyle{plain}
\newtheorem{thm}{Theorem}[section]
\newtheorem{lem}[thm]{Lemma}
\newtheorem{pro}[thm]{Proposition}
\newtheorem{cor}[thm]{Corollary}
\theoremstyle{definition}

\newtheorem{rem}[thm]{Remark}
\newtheorem{exa}[thm]{Example}

\begin{document}
\title{On a conjecture of Wilf}
\author{Stefan De Wannemacker, Thomas Laffey, and Robert Osburn}

\address{Faculty of Applied Economic Sciences, University of Antwerp, Prinsstraat 13, 2000 Antwerp, Belgium}
\email{stefan.dewannemacker@ua.ac.be}

\address{School of Mathematical Sciences, University College Dublin, Belfield, Dublin 4, Ireland}

\email{thomas.laffey@ucd.ie}
\email{robert.osburn@ucd.ie}


\date{January 26, 2007}
\maketitle

\begin{abstract}
Let $n$ and $k$ be natural numbers and let $S(n,k)$ denote the Stirling numbers of the second kind. It is a conjecture of Wilf that the alternating sum

\begin{center}
$\displaystyle \sum_{j=0}^{n} (-1)^{j} S(n,j)$
\end{center}

\noindent is nonzero for all $n>2$. We prove this conjecture for all $n\not \equiv 2$ and $\not \equiv 2944838 \bmod 3145728$ and discuss applications of this result to graph theory, multiplicative partition functions, and
the irrationality of $p$-adic series.

\end{abstract}

 \section{Introduction}

Let $n$ and $k$ be natural numbers. The Stirling numbers $S(n,k)$ of the second kind are given by
\[x^n=\sum_{k=0}^{\infty}S(n,k)(x)_k,\] where
$(x)_k:=x(x-1)(x-2)\ldots(x-k+1)$ for $k \in \mathbb{N}\setminus
\{0\}$ and $(x)_0:=1.$ $S(n,k)$ is the number of ways in
which it is possible to partition a set with $n$ elements into
exactly  $k$\ nonempty subsets. Consider the alternating sum

$$f(n) := \sum_{j=0}^n (-1)^{j}S(n,j).$$\\

\noindent The first few terms in the sequence of integers $\{f(n)\}_{n\geq 0}$ are as follows:

\begin{center}
$1$, $-1$, $0$, $1$, $1$, $-2$, $-9$, $-9$, $50$, $267$, $413$, $-2180$, $-17731$, $-50533$, $110176$, $\dotsc$
\end{center}

\noindent This is sequence $A000587$ of Sloane \cite{sloane}. This
sequence appears in Example 5(ii), Section 8, Chapter 3 in Ramanujan's second
notebook (see page 53 of \cite{berndt}) and has been subsequently investigated by Beard \cite{beard}, Harris
and Subbarao \cite{hs}, Uppuluri and Carpenter \cite{uc},
Kolokolnikova \cite{k}, Layman and Prather \cite{lp}, Subbarao and
Verma \cite{sv}, Yang \cite{yang}, Klazar \cite{klazar1}, and Murty
and Sumner \cite{ms}.

Wilf has conjectured (see \cite{klazar}) that $f(n) \neq 0$ for all
$n>2$. So, the only value of $n$ for which $f(n)$ vanishes would be
$n=2$. The best known result in this direction is that of Yang
\cite{yang}. In \cite{yang}, the author adapted an approach of de
Bruijn \cite{bruijn} concerning the saddle point method and used
exponential sum estimates from \cite{kn} to show that the number of
$n \leq x$ with $f(n)=0$ is $O(x^{2/3})$ where the implied constant is not explicitly computed.
Recently, Murty and Sumner have taken a different approach in proving the non-vanishing of $f(n)$. In \cite{ms}, the authors use the congruence

\begin{center}
$\displaystyle f(n) \equiv \sum_{j=0}^{n} S(n,j) \equiv B_{n} \bmod 2$,
\end{center}

\noindent properties of the Bell numbers $B_{n}$, and of $\zeta_{3}$, a cube root of unity, to prove the following result.

\begin{thm} \label{ms} If $n \not \equiv 2 \bmod 3$, then $f(n) \neq 0$.
\end{thm}

The purpose of this paper is to extend Theorem \ref{ms} as follows.

\begin{thm} \label{main} If $n \not \equiv 2$ and $\not \equiv 2944838 \bmod 3145728$, then $f(n) \neq 0$.
\end{thm}

The paper is organized as follows. In Section 2, we use generating functions and properties of finite fields to prove a general congruence for $f(n)$. This congruence (see Proposition \ref{maincon}) combined with computer calculations (see the Appendix) yields a proof of Theorem 1.2. In Section 3, we give a brief discussion of other congruences for $f(n)$. In particular we prove a general congruence for $f(n)$ modulo $p$ where $p$ is a prime (see Proposition \ref{othercon}). This result generalizes the congruences given by Lemmas 9 and 10 in \cite{ms}. We conclude Section 3 by mentioning another approach to the general congruence for $f(n)$ using a certain set of recursively defined polynomials. We relate these polynomials to $f(n)$ and use this relationship to give an alternative proof of Proposition \ref{maincon}. In Section 4, we discuss how Theorem \ref{main} has applications in three distinct areas of mathematics, namely graph theory, multiplicative partition functions, and to the irrationality of $p$-adic series.

\section{Proof of Theorem 1.2}

The proof of Theorem \ref{main} contains two key steps. We first derive the generating function for $f(n)$, then use this expression to determine when $f(n)$ has a period of $N$ modulo $m$ where $N$ and $m$ are positive integers. We first require the following well-known
property of Stirling numbers of the second kind, namely (see page 34 in \cite{stan})

\begin{equation}\label{gfs}
\sum_{n\geq k} S(n,k) x^{n} = \frac{x^k}{(1-x)(1-2x) \dotsm (1-kx)}.
\end{equation}

\noindent For more details and basic results on Stirling numbers of the second kind we refer
the reader to \cite{COM}, \cite{GKP}, or \cite{stan}. Recent applications of
$S(n,k)$ include computing annihilating polynomials for quadratic
forms \cite{s1}. Further information on these applications can be
found in \cite{s3}. Using (\ref{gfs}), we derive an expression for the generating function of $f(n)$.

\begin{lem} \label{gff} For any positive integer $m$, the generating function $F(x)$ of $f(n)$ is the following rational function modulo $m$.

\begin{equation} \label{fpd}
F(x):=\sum_{n\geq 0} f(n) x^n \equiv \frac{Q(x)}{(1-x)(1-2x) \dotsm (1-(m-1)x) - (-1)^{m} x^m} \bmod m,
\end{equation}

\noindent where $Q(x)$ is a polynomial modulo $m$ given by

$$
Q(x):= \Biggl (\sum_{k=0}^{m-1} \frac{(-1)^k x^k}{(1-x)(1-2x) \dotsm (1-kx)} \Biggr) (1-x)(1-2x)\dotsm(1-(m-1)x).
$$
\end{lem}

\begin{proof}
We begin by multiplying both sides of (\ref{gfs}) by $(-1)^k$ and summing over $k$ to obtain

\begin{equation} \label{genf}
F(x)=\sum_{n\geq 0} f(n) x^n = \sum_{k \geq 0} \frac{(-1)^k x^k}{(1-x)(1-2x) \dotsm (1-kx)}.
\end{equation}

\noindent Now computing $F(x)$ modulo $m$ yields

$$
\begin{aligned}
& F(x) \\
&\equiv  \Biggl (\sum_{k=0}^{m-1} \frac{(-1)^k x^k}{(1-x)(1-2x) \dotsm (1-kx)} \Biggr) \cdot  \Biggl (\sum_{i=0}^{\infty} \Bigl( \frac{(-1)^m x^m}{(1-x)(1-2x) \dotsm (1-(m-1)x)}\Bigr)^{i} \Biggr) \bmod m \\
& \equiv \Biggl (\sum_{k=0}^{m-1} \frac{(-1)^k x^k}{(1-x)(1-2x) \dotsm (1-kx)} \Biggr) \cdot  \Biggl (1 - \frac{(-1)^m x^m}{(1-x)(1-2x) \dotsm (1-(m-1)x)}\Biggr)^{-1} \bmod m \\
& \equiv \frac{Q(x)}{(1-x)(1-2x) \dotsm (1-(m-1)x) - (-1)^{m} x^m} \bmod m, \\
\end{aligned}
$$

\noindent where $Q(x)$ is defined as above.
\end{proof}

\begin{rem} \label{dx}
Given a positive integer $m$, we now explain one way to compute a period $N$ for $f(n)$ modulo $m$. Consider

$$D(x):=(1-x)(1-2x) \dotsm (1-(m-1)x) - (-1)^{m} x^m,$$

\noindent which is the denominator of $F(x)$ via Lemma \ref{gff}. Note that $F(x)$ is proper, i.e., the degree of $Q(x)$ is less than the degree of $D(x)$. Let $\alpha$ be a root of $D(x)$ modulo $m$ and view $\alpha$ as the representative of $x$ in the ring $\mathbb{Z}_{m}[x]/\langle D(x) \rangle$. Let

$$
\Gamma(\alpha):=\{ a_{0} + a_{1} \alpha + \cdots a_{m-1} {\alpha}^{m-1} : a_{i} \in \mathbb{Z}_{m} \}.
$$

\noindent Then $\Gamma(\alpha)$ forms a finite semi-group under multiplication. Define $\Gamma^{*}(\alpha)$ to be the set of invertible elements in $\Gamma(\alpha)$. Then $\Gamma^{*}(\alpha)$ forms a finite group. Moreover, let $g(x)=(1-D(x))/x$ and note that $g(\alpha)$ is a polynomial in $\alpha$ of degree at most $m-1$ and hence belongs to $\Gamma(\alpha)$. Also we have $g(\alpha)\alpha =1$
and so $\alpha$ belongs to $\Gamma^{*}(\alpha)$. As the order of $\alpha$ divides $|\Gamma^{*}(\alpha)|$,

$$
\displaystyle \alpha^{|\Gamma^{*}(\alpha)|} =1,
$$

\noindent and so $\alpha$ is a root of $x^{|\Gamma^{*}(\alpha)|} - 1$. Since this is true for all roots of $D(x)$, we get

$$
1-x^{N} \equiv D(x)M(x) \bmod m
$$

\noindent where $M(x) \in \mathbb{Z}_{m}[x]$ and $N$ is the least common multiple of the $|\Gamma^{*}(\alpha)|$ as $\alpha$ ranges over the roots of $D(x)$. Now if $1-x^{N} \equiv D(x)M(x) \bmod m$, then observe that the proper rational function

$$
F(x) \equiv \frac{Q(x)M(x)}{1-x^N} \bmod m
$$

\noindent has a period of $N$ upon multiplying both sides by $1-x^{N}$ and comparing coefficients. This in turn implies that $f(n)$ has a period of $N$ modulo $m$. For example, if $m=2$, then

$$
F(x) \equiv \frac{1}{x^2 + x + 1} \bmod 2.
$$

\noindent By multiplying both the numerator and denominator by $x+1$, we obtain

$$
F(x) \equiv \frac{x+1}{x^3 -1} \bmod 2.
$$

\noindent Thus $f(n) \equiv f(n+3) \bmod 2$ and so we recover Theorem 1.1 as $f(1)$ and $f(3)$ are odd.

\end{rem}

We are now in a position to prove a general congruence for $f(n)$.

\begin{pro} \label{maincon} Let $n$, $h \in \mathbb{N}$. Then
$$f(n)\equiv f(n+3\cdot4^{h-1}) \bmod 2^{h}.$$
\end{pro}

\begin{proof} We first work over the field $\mathbb{F}_{2}$. Let $m=2^h$ with $h\geq 1$. By Remark \ref{dx}, it is sufficient to find a positive integer $N$ such that $\alpha^{N} \equiv 1 \bmod m$ whenever $D(\alpha) \equiv 0 \bmod m$. Let $\alpha$ be a root of $D(x)$. Then

$$
\begin{aligned}
D(\alpha) &= (1-\alpha)(1-2\alpha)\dots(1-(2^{h} - 1)\alpha) - (-1)^{2^{h}} \alpha^{2^{h}} \\
& \equiv (1-\alpha)^{2^{h-1}} + \alpha^{2^h} \bmod 2 \\
& \equiv 1 + \alpha^{2^{h-1}} + \alpha^{2^{h}} \bmod 2\\
& \equiv 0 \bmod 2
\end{aligned}
$$

\noindent and thus $\alpha^{3\cdot 2^{h-1}} \equiv 1 \bmod 2$. So we have $\alpha^{3\cdot 2^{h-1}} \equiv 1 + 2w \bmod 2^h$ for some $w \in \mathbb{Z}$. Then

$$
\alpha^{3 \cdot 2^{h-1} \cdot 2^{h-1}} \equiv (1+2w)^{2^{h-1}}
\equiv 1+ 2^{h} w + \binom{2^{h-1}}{2} 2^2 w^2 + \dotsb +
(2w)^{2^{h-1}} \bmod 2^h.
$$

\noindent As $\binom{2^{h-1}}{t} \equiv 0 \bmod 2^{h}$ for all $1 \leq t \leq 2^{h-1}$, we deduce

$$
\alpha^{3\cdot 4^{h-1}} \equiv 1 \bmod 2^h
$$

\noindent and thus $3\cdot 4^{h-1}$ is a period for $f(n)$ modulo $2^h$.

\end{proof}

We can now prove Theorem \ref{main}

\begin{proof}
For every fixed value of $h \geq 1$ one can use Proposition
\ref{maincon} to compute the values of $n$ in the interval $[0$,
$3\cdot4^{h-1}-1 ]$ for which the 2-adic valuation of $f(n)$ is at
least $h$. These values will yield the only possible cases $\bmod
~3\cdot4^{h-1}$ for which Wilf's conjecture might fail, the
so-called ``open'' cases. For large values of $h$, the computer
program given in the Appendix can be used for this purpose. In
particular, take $h=22$ and consider the set

\begin{center}
 $N_{22} := \{l \in \mathbb{N}\ :   l < 3\cdot4^{21}$ and $f(l) \not \equiv 0
\bmod 2^{22} \}$.
\end{center}

\noindent The congruence
\[f(n)\equiv f(n+3\cdot4^{21}) \bmod 2^{22} \\
\]
implies that

\[f(N) \neq 0 \ \]

\noindent for all $N \equiv l \bmod 3\cdot 4^{21}$ where $l \in
N_{22} $. In particular, since $f(n)\equiv 0 \bmod 2^{22} $ only for
the values $n \equiv 2$ and $\equiv 2944838 \bmod 3145728$ when
$n<3\cdot 4^{21}$, this implies that if $n\not\equiv 2$ and
$\not\equiv 2944838 \bmod 3145728$, then $f(n) \neq 0$ and the
result follows.
\end{proof}

In the table below we have listed the ``open'' cases for values of $h\leq
22$.

\begin{center}
\begin{tabular}{|c|c|c|}
\multicolumn{3}{r}{} \\ \hline
$h$ & Open cases & $\mod$ \\
\hline
1& $2 $& $3$\\
2& $2, 11 $&$ 12$\\
3& $2 $& $ 12$\\
4& $2 $& $ 12$\\
5& $2 $& $ 12$\\
6& $2, 38 $& $ 48$\\
7& $2, 38 $& $ 96$\\
8& $2, 134 $& $ 192$\\
9& $2, 326 $& $ 384$\\
10& $2, 326 $& $ 768$\\
11& $2, 326 $& $ 1536$\\
12& $2, 1862 $& $ 3072$\\
13& $2, 1862 $& $ 6144$\\
14& $2, 8006 $& $ 12288$\\
15& $2, 20294 $& $ 24576$\\
16& $2, 44870 $& $ 49152$\\
17& $2, 94022 $& $ 98304$\\
18& $2, 192326 $& $ 196608$\\
19& $2, 192326 $& $ 393216$\\
20& $2, 585542 $& $ 786432$\\
21& $2, 1371974 $& $ 1572864$\\
22& $2, 2944838 $& $ 3145728$\\
\hline
\end{tabular}
\end{center}

\section{Other congruences}

The purpose of this section is two-fold. We first discuss how Remark \ref{dx} can also be used to prove other interesting congruences for $f(n)$. Secondly, we provide an alternative approach to proving congruences for $f(n)$ using a recursively defined set of polynomials. We begin with an immediate application of Remark \ref{dx}.

\begin{pro} \label{othercon} Let $n$, $h \in \mathbb{N}$ and $p$ be an odd prime. Then

\begin{equation} \label{modp}
\displaystyle f(n) \equiv f(n + 2\tfrac{p^{p} - 1}{p-1}) \bmod p
\end{equation}

\noindent and

\begin{equation} \label{modph}
f(n) \equiv f(n + \tfrac{2p^{2h-2}(p^{p} -1)}{(p-1)}) \bmod p^{h}.
\end{equation}

\end{pro}

\begin{proof} We work over the field $\mathbb{F}_{p}$. By Fermat's Little Theorem for finite fields, the denominator $D(x)$ can be simplified, namely

$$
D(x) = (1-x)(1-2x) \dotsm (1-(p-1)x) + x^{p} \equiv 1- x^{p-1} + x^{p} \bmod p.
$$

\noindent By Remark \ref{dx}, we assume that $\alpha$ is a root of $D(x)$ and let $\beta=1/\alpha$. Note that the period of $\alpha$ is the same as the period of $\beta$. One can then check that

\begin{equation} \label{beta}
\beta(\beta -1)(\beta - 2) \dotsm (\beta - p + 1) + 1 \equiv \beta^{p} - \beta + 1 \equiv 0 \bmod p.
\end{equation}

\noindent We now show by induction on $i$ that $\beta^{p^{i}} \equiv \beta - i \bmod p$. The result holds for $i=0$ and $i=1$ by (\ref{beta}). Assume $\beta^{p^{i}} \equiv \beta - i \bmod p$. Then

$$
\beta^{p^{i+1}} \equiv \Big( \beta^{p^{i}} \Big)^{p} \equiv (\beta - i)^{p} \equiv \beta^{p} - i \equiv \beta - (i+1) \bmod p.
$$

\noindent This proves the claim. Now applying this claim and (\ref{beta}), we have

$$
\beta^{1+ p + p^2 + \dotsm + p^{p-1}} \equiv \beta(\beta -1)(\beta - 2) \dotsm (\beta - p + 1) \equiv -1 \bmod p.
$$

\noindent Therefore $\frac{2(p^{p} - 1)}{p-1}$ is a period of $\beta$ and hence is a period of $\alpha$. By Remark \ref{dx}, (\ref{modp}) then follows. The proof of (\ref{modph}) is similar to that of Proposition \ref{maincon} and is left to the reader.

\end{proof}

\begin{rem}
One can ask for the minimal periods for $f(n)$ modulo $m$. In the table below, we compute the minimal periods for $f(n)$ modulo $m$
for small values of $m$. The values in this table follow from
Propositions \ref{maincon}, \ref{othercon},
and numerical work. Note that $f(n)$ does not have minimal period $3\cdot 4^{h-1}$ modulo $2^{h}$ as can be seen for $h=3$ (see Remark \ref{notmin}). In general, we conjecture that the minimal period for $f(n)$ modulo $p^{h}$ where $h \geq 1$ is the one given by (\ref{modph}). We would like to point out (thanks to the referee) that the congruences for $f(n)$ are completely analogous to congruences for the Bell numbers. In particular, it is well known that for prime $p$, the Bell numbers are periodic with minimal period dividing $\frac{p^{p} - 1}{p-1}$ and that this seems to be the minimal period. No one has been able to prove this claim. For further information regarding congruences for Bell numbers, please see \cite{car}, \cite{lps}, and \cite{wag}. \\

\begin{center}
\begin{tabular}{|l|r||l|r|}
\hline
$m$& Minimal period&$m$& Minimal period\\
\hline
&&&\\
2&3&10&398310\\
&&&\\
3& $3^3-1$&11& $\frac{11^{11}-1}{5}$\\
&&&\\
4& $3\cdot 4$&12&1560\\
&&&\\
5& $\frac{5^5-1}{2}$&13&$\frac{13^{13}-1}{6}$\\
&&&\\
6&390&14&17294382\\
&&&\\
7& $\frac{7^7-1}{3}$&15&81091300290\\
&&&\\
8& $\frac{3\cdot 4^2}{2}$&16&$3 \cdot 4^3$\\
&&&\\
9&$\frac{2\cdot 3^{2} (3^{3} -1)}{3-1}$&&\\
&&&\\
\hline
\end{tabular}
\end{center}

\end{rem}

We now turn to an alternative approach to proving Proposition \ref{maincon}. Consider the set of polynomials defined in the following recursive
way:

\begin{align*}
P_0(X) := &\ 1\\
P_n(X) := &\ XP_{n-1}(X) - P_{n-1}(X+1),\ n\geq 1.\\
\end{align*}

\begin{exa}
\begin{eqnarray*}
P_1(X) &=& X-1,\\
P_2(X) &=& X^2-2X,\\
P_3(X) &=& X^3-3X^2+1,\\
P_4(X) &=& X^4-4X^3+4X+1,\\
P_5(X) &=& X^5 -5X^4+10X^2+5X-2.\\
\end{eqnarray*}
\end{exa}

\noindent The generating function of the $P_{n}$'s is given by

\begin{equation} \label{gfp}
P(X,t):= \sum_{n \geq 0} P_{n}(X) t^n = \sum_{j\geq 0}  \frac{(-1)^j t^j}{(1-Xt)(1-(X-1)t) \dotsm (1-(X-j)t)}.
\end{equation}

\noindent To see this, multiply the recurrence for the $P_{n}$'s by $t^n$ and sum over $n$ to get the functional equation

$$
P(X,t) = \frac{1-tP(X+1, t)}{1-Xt},
$$

\noindent which is satisfied by (\ref{gfp}). We now relate these polynomials to $f(n)$ and prove a recursive formula. Precisely, we have

\begin{pro} \label{prop2}
Let $n \in \mathbb{N}$. Then \\
$(i)$ $\displaystyle f(n)=P_{n}(0).$\\
$(ii)$ $\displaystyle P_{n}(X) = \sum_{j=0}^{n} \binom{n}{j} f(n-j) X^{j}.$ \\
$(iii)$ $\displaystyle -f(n+1)= \sum_{j=0}^{n} \binom{n}{j} f(n-j).$
\end{pro}

\begin{proof} Taking $X=0$ in (\ref{gfp}) and using (\ref{genf}) yields $(i)$. Now $(ii)$ follows from comparing the coefficient of $t^n$ in (\ref{gfp}) and using (\ref{gfs}). Finally, by $(ii)$, we have

$$
P_{n}(1) = \sum_{j=0}^{n} \binom{n}{j} f(n-j).
$$

\noindent Then $(iii)$ follows since $f(n+1) = P_{n+1}(0) = - P_n(1)$. We note that observation $(iii)$ was originally made in the context of multiplicative partition functions (see \cite{sv}).

\end{proof}

\begin{rem} It has been numerically verified that $P_{n}(X)$ is irreducible over $\mathbb{Z}$ for all $5 < n \leq 200$. We believe that $P_{n}(X)$ is irreducible over $\mathbb{Z}$ for all $n>5$. It is not immediately clear that the methods of \cite{c}, \cite{ft}, or \cite{es} can be suitably adapted to prove this claim. Note that this claim implies Wilf's conjecture as the constant term of $P_{n}(X)$ is $f(n)$.
\end{rem}

We now prove the following useful properties of the polynomials $P_n(X)$.

\begin{pro} \label{stuff}
Let $k$ be a positive integer. Let
\begin{align*}
f_{k}(X,Y)&:=(X-Y)(X+1-Y)\cdots(X+k-1-Y)\\
&=\sum_{r=0}^{k}a_{r,k}(X)Y^r
\end{align*}
where $a_{r,k}(X) \in \mathbb{Z}[X].$ Then for all $n \in \mathbb{N}$,
\[P_{n}(X+k)=\sum_{r=0}^{k}a_{r,k}(X)P_{n+r}(X).\]
\end{pro}

\begin{proof}
We proceed by induction on $k$. When $k=1$, the result states
$P_n(X+1)=XP_n(X)-P_{n+1}(X)$ and this is the recurrence
relation for the polynomials $P_{n}(X).$ Assume the result holds for $k$. Then
\begin{align*}
P_n(X+k+1) &= \sum_{r=0}^{k}a_{r,k}(X+1)P_{n+r}(X+1)\\
&=\sum_{r=0}^{k}a_{r,k}(X+1)\left(XP_{n+r}(X)-P_{n+r+1}(X)\right)\\
&=\sum_{r=0}^{k}\left(Xa_{r,k}(X+1)P_{n+r}(X)-a_{r,k}(X+1)P_{n+r+1}(X)\right).\\
\end{align*}
For $0\leq t\leq k$, the coefficient of $P_{n+t}(X)$ is
\[Xa_{t,k}(X+1)-a_{t-1,k}(X+1).\]

\noindent Thus

\begin{align*}
f_{k+1}(X,Y)&=(X-Y)f_{k}(X+1,Y)\\
&=(X-Y)\sum_{r=0}^{k}a_{r,k}(X+1)Y^r\\
&=\sum_{r=0}^{k}Xa_{r,k}(X+1)Y^r-\sum_{r=1}^{k+1}Xa_{r-1,k}(X+1)Y^{r}.\\
\end{align*}
So $a_{t,k+1}(X)=Xa_{t,k}(X+1)-a_{t-1,k}(X+1)$ and
\[
P_n(X+k+1)\ =\ \sum_{r=0}^{k+1}a_{r,k+1}(X)P_{n+r}(X).
\]
\end{proof}

\begin{cor}\label{modk_cor}
Let $n$, $k \in \mathbb{N}$. Then
\[f(n)\equiv \sum_{r=1}^{k}a_{r,k}(0)f(n+r) \mod k \]
where
\[(X-Y)(X+1-Y)\cdots(X+k-1-Y)=\sum_{r=0}^{k}a_{r,k}(X)Y^r.\]
\end{cor}

\begin{proof}
From Proposition \ref{stuff} and $a_{0,k}(0)=0$, it follows that
\begin{align*}
P_{n}(k)=&\ \sum_{r=1}^{k}a_{r,k}(0)P_{n+r}(0)\\
=&\ \sum_{r=1}^{k}a_{r,k}(0)f(n+r).\\
\end{align*}

\noindent The result now follows from part $(i)$ of Proposition \ref{prop2} and the fact that

\[ P_{n}(k) \equiv P_n(0) \mod k.\]
\end{proof}

We can now give an alternative proof of Proposition \ref{maincon}.

\begin{proof}
Corollary \ref{modk_cor} for $k=2^{h}$ gives
\[f(n) \equiv \sum_{r=1}^{2^{h}}a_{r,2^{h}}(0) f(n+r) \mod 2^{h}\]
and, in particular,
\[f(n+2^{h}) \equiv f(n) -\sum_{r=1}^{2^{h} - 1}a_{r,2^{h}}(0) f(n+r) \mod 2^{h}.\]
So we have
\[\begin{pmatrix}
f(n+2^{h})\\
f(n+2^{h}-1)\\
f(n+2^{h}-2)\\
\vdots\\
f(n+1)\\
\end{pmatrix}
\equiv A\begin{pmatrix}
f(n+2^{h}-1)\\
f(n+2^{h}-2)\\
f(n+2^{h}-3)\\
\vdots\\
f(n)\\
\end{pmatrix}
\mod 2^{h}\] where \[A=\begin{pmatrix}
-a_{2^h-1,2^h}(0) & -a_{2^h-2,2^h}(0) & -a_{2^h-3,2^h}(0)& \cdots &-a_{1,2^{h}}(0)& 1 \\
1   & 0   & 0  & \cdots &0&0\\
0   & 1   & 0  & \cdots &0&0\\
\vdots& \vdots & \vdots & & \vdots & \vdots \\
0   & 0   & 0  & \cdots &1&0\\
\end{pmatrix}
.\]

\noindent Note that $A$ is the companion matrix of the polynomial

\begin{center}
$c(Y)=Y(Y-1)\cdots(Y-2^{h+1}+1)-1.$
\end{center}

\noindent Now

\begin{center}
$c(Y)\equiv (Y(Y+1))^{2^{h-1}}+1 \equiv (Y^2+Y+1)^{2^{h-1}} \mod 2$.
\end{center}

\noindent Over
$\mathbb{F}_{2}$, $A$ is non-derogatory (see \cite{Brown}, 7.20) and
has Jacobson canonical form (see
\cite{Jacobson}, page 72)
\[J=\begin{pmatrix}X & N & & &\\
& X & N &&\\
&&\ddots&\ddots&\\
&&&&N\\
&&&&X\\
\end{pmatrix}\]
where $X=\begin{pmatrix}1&1\\1&0\\\end{pmatrix}$ is the companion
matrix of $Y^2+Y+1$ and $N=\begin{pmatrix}0&0\\1&0\\\end{pmatrix}.$
Let $I_{s}$ be the $s \times s$ identity matrix where $s \geq 1$.
Some calculation shows that
\[c(J)=\begin{pmatrix}0 & I_{2} & & &\\
& 0 & I_{2} &&\\
&&\ddots&\ddots&\\
&&&0&I_{2} \\
&&&&0\\
\end{pmatrix}\] and
\[Z:=J^3-I_{2^{h}}=(J-I_{2^{h}})c(J)=\begin{pmatrix}0 & X^2 & N& &&\\
& 0 & X^2 &N&&\\
&&\ddots&\ddots&\ddots&\\
&&&0&X^2&N\\
&&&&0&X^2\\
&&&&&0\\
\end{pmatrix}.\]
\noindent The matrix $Z$ has the property that

\begin{center}
$Z^{2^{h-1}} = 0$.
\end{center}

\noindent So over $\mathbb{F}_{2}$,

\begin{align*}
J^{3\cdot2^{h-1}} & = (I_{2^{h}} + Z)^{2^{h-1}}\\
& = I_{2^{h}} + 2^{h-1}Z + \ldots + Z^{2^{h-1}}\\
& = I_{2^{h}}.\\
\end{align*}
Hence $A^{3\cdot2^{h-1}}$ is similar to a matrix of the form  $I_{2^{h}}+2W$ for
a matrix $W$ over $\mathbb{F}_{2^{h}}$. So
\begin{align*}
A^{3\cdot2^{h-1}\cdot 2^{h-1}}&=(I_{2^{h}}+2W)^{2^{h-1}}\\
&= I_{2^{h}} + 2^hW+ {2^{h-1}\choose 2}4W^2+\ldots+(2W)^{2^{h-1}}.\\
\end{align*}
Since $\displaystyle {2^{h-1}\choose t}2^t\equiv 0 \mod 2^h$, for all $1\leq t
\leq 2^{h-1}$, we have
\[A^{3\cdot2^{2h-2}}\equiv I_{2^{h}} \mod 2^h.\]
In other words,  we have
\begin{align*}
\begin{pmatrix}
f(n+2^{h})\\
f(n+2^{h}-1)\\
f(n+2^{h}-2)\\
\vdots\\
f(n+1)\\
\end{pmatrix}
&\equiv A\begin{pmatrix}
f(n+2^{h}-1)\\
f(n+2^{h}-2)\\
f(n+2^{h}-3)\\
\vdots\\
f(n)\\
\end{pmatrix} \bmod 2^{h} \\
\\
&\equiv A^2\begin{pmatrix}
f(n+2^{h}-2)\\
f(n+2^{h}-3)\\
f(n+2^{h}-4)\\
\vdots\\
f(n-1)  \\
\end{pmatrix}  \bmod 2^{h} \\
\\
&\vdots\\
\\
&\equiv \begin{pmatrix}
f(n+2^{h}-3\cdot2^{2h-2})\\
f(n+2^{h}-3\cdot2^{2h-2}-1)\\
f(n+2^{h}-3\cdot2^{2h-2}-2)\\
\vdots\\
f(n-3\cdot2^{2h-2}-2^{h})\\
\end{pmatrix} \bmod 2^h .
\end{align*} \\
Comparing the elements in the matrix yields the desired congruence.
\end{proof}

\begin{rem} \label{notmin}
Using this polynomial approach, one can show for instance that $3\cdot 4^{h-1}$ is not a minimal period modulo $2^{h}$, $h \geq 1$, for $f(n)$. Namely, one can use the definition of the $P_{n}$'s and Proposition \ref{prop2} to check that

$$
P_{n+48}(0) \equiv P_{n+24}(0) \bmod 8
$$

\noindent or equivalently

$$
f(n) \equiv f(n+24) \bmod 8.
$$

\end{rem}

\section{Applications}
In this section we consider applications of Theorem \ref{main} to graph theory, multiplicative partition functions, and to the irrationality of a $p$-adic series.

\subsection{Graph Theory}

A {\it simple graph} $G$ consists of a non-empty finite set $V(G)$ of {\it vertices} and a finite set $E(G)$ of distinct unordered pairs of distinct elements of $V(G)$ called {\it edges}. We say that two vertices $v$, $w \in V(G)$ are {\it adjacent} if there is an edge $(v,w) \in E(G)$ joining them. A graph for which $E(G)$ is empty is called the {\it null graph} and is denoted by $N_{n}$ where $n$ is the number of vertices. A {\it complete graph} is a simple graph in which each pair of distinct vertices are adjacent. The complete graph on $n$ vertices is denoted by $K_{n}$. If the vertex set of a graph $G$ can be partitioned into two disjoint sets $A$ and $B$ so that each edge of $G$ joins a vertex of $A$ and a vertex of $B$, then $G$ is called a {\it bipartite graph}. A {\it complete bipartite graph} is a bipartite graph in which each vertex of $A$ is joined to each vertex of $B$ by just one edge. The complete bipartite graphs are denoted by $K_{r,s}$ where $r$ and $s$ are the cardinalities of $A$ and $B$ respectively.

Let $G$ be a simple graph with $n$ vertices. One can associate to $G$ many polynomials whose properties yield structure theorems of isomorphism classes of graphs. In the vast literature, one can study, for example, the {\it chromatic polynomial}, {\it Tutte polynomial}, {\it interlace polynomials}, {\it cover polynomials} of digraphs, and the {\it matching polynomial} of a graph. In this section we take a closer look at the matching polynomial of certain bipartite graphs.

A {\it $k$-matching} in a graph $G$ is a set of $k$ edges, no two of which have a vertex in common. We denote the number of $k$-matchings in $G$ by $p(G,k)$. We set $p(G,0)=1$ and define the {\it matching polynomial} of $G$ by

\begin{center}
$\displaystyle \mu(G,X):= \sum_{k \geq 0} (-1)^{k} p(G,k) X^{n-2k}$.
\end{center}

\noindent Some examples of matchings polynomials are

\begin{center}
$\mu(N_{n},X)=X^n$,
\end{center}

\begin{center}
$\mu(K_{n},X)=\displaystyle \sum_{k \geq 0} (-1)^{k} \frac{n!}{k!(n-2k)! 2^{k}} X^{n-2k}$,
\end{center}

\noindent and

\begin{center}
$\mu(K_{n,n},X)=\displaystyle \sum_{k \geq 0} (-1)^{k} {\binom{n}{k}}^2 k!  X^{n-2k}$.
\end{center}

\noindent The study of matching polynomials has been a focus of research over the last twenty five years. For further details regarding properties of matching polynomials, the reader should consult \cite{bf}, \cite{dg}, \cite{farrell}, \cite{g}, \cite{godsil}, \cite{gg}, \cite{gg1}, or \cite{l}. As we are interested in the roots of $\mu(G,X)$,  we recall some general results.

\begin{pro} Let $G$ be a graph with $n$ vertices. Then \\
(i) The zeros of $\mu(G,X)$ are real. \\
(ii) The zeros of $\mu(G,X)$ are symmetrically distributed about the origin. \\
\end{pro}

\begin{proof}
For (i), see Corollary 1.2 or Lemma 4.3 in \cite{godsil}. If $n$ is even, then $\mu(G,X)$ can be written as a polynomial in $X^2$. If $n$ is odd, then $X^{-1} \mu(G,X)$ can be expressed as a polynomial in $X^2$. Thus (ii) follows.
\end{proof}

\noindent Further results on roots of matching polynomials can be found in \cite{fr}, \cite{gwalks}, \cite{gg}, \cite{gg1}, \cite{gut}, or \cite{hl}. For our purposes, we consider the following bipartite graph. Let $T(n)$ be the graph with vertex set $\{1,\dotsc, n\} \cup \{1^{\prime}, \dotsc, n^{\prime} \},$ where $i$ is adjacent to $j^{\prime}$ if and only if $i>j$.  Thus $T(n)$ has $2n$ vertices. For $n=3$, one can check that $p(T(3),1)=3$, $p(T(3),2)=1$, $p(T(3),3)=0$, and thus

$$\displaystyle \mu(T(3), X)=X^2(X^2 - X -1)(X^2 + X -1).$$

We now relate the matching polynomial of $T(n)$ to Stirling numbers of the second kind $S(n,k)$.

\begin{pro} For the graph $T(n)$, we have
\begin{center}
$\displaystyle \mu(T(n),X)=\sum_{k=0}^{n} (-1)^k S(n,n-k) X^{2n-2k}$.
\end{center}
\end{pro}

\begin{proof}
We briefly sketch the proof as given in \cite{godsil}. For another proof, see the solution to Problem 4.31 in \cite{lovasz}. The idea is to consider a bijection from the set of $k$-matchings of $T(n)$ to a certain set of directed graphs. Thus counting the number of such directed graphs yields $p(T(n),k)$ and thus $\mu(T(n),X)$. Each matching in $T(n)$ determines a directed graph with vertex set $N=\{1, \dotsc, n\}$ with arc $(i,j)$ for each edge $\{i,j^{\prime}\}$ in the matching and a loop on each vertex $j$ not in the matching. Now, each vertex component is a directed path with a loop on its last vertex. As there is an arc from $i$ to $j$ in the directed graph only if $i\geq j$, the graph is determined by the vertex set of each component. Thus the number of such directed graphs with $c$ components is $S(n,c)$. Note that $c$ equals the number of loops and decreases by $1$ for each edge in the original matching. Hence $c=n-k$ where $k$ equals the number of edges in the matching.
\end{proof}

\noindent Each of the polynomials $\mu(T(n),X)$ contains $X^2$ as a factor and thus is reducible. We thus consider the roots of the polynomial $\displaystyle \frac{1}{X^2} \mu(T(n), X)$. This corresponds to removing
the vertices $1$ and $n^{\prime}$ in the graph $T(n)$. As a result of Theorem \ref{main}, we immediately have

\begin{cor} For $n \not \equiv 2$ and $\not\equiv 2944838 \bmod 3145728$, $1$ is not a root of $\displaystyle \frac{1}{X^2} \mu(T(n), X)$.
\end{cor}

\begin{rem}
We conjecture that $1$ is not a root of $\displaystyle \frac{1}{X^2} \mu(T(n), X)$ for  $n \equiv 2$ and $\equiv 2944838 \bmod 3145728$, and, more generally, that $\displaystyle \frac{1}{X^2} \mu(T(n), X)$ is irreducible over $\mathbb{Z}$ for every $n>3$. This last statement has been numerically verified for all $3 < n \leq 500$. Note that this statement implies Wilf's conjecture.
\end{rem}

\subsection{Multiplicative partition functions}

Multiplicative partition functions count the number of representations of a given positive integer $m$ as a product of positive integers. For a well-written survey of techniques for enumerating product representations, please see \cite{km}. Suppose the canonical prime factorization of $m$ is given by

\begin{center}
$\displaystyle m=p_{1}^{r_{1}} \dotsc p_{n}^{r_{n}}.$
\end{center}

\noindent The succession of integers $r_{1}$, $r_{2}$, $\dotsc r_{n}$, when arranged in descending order of magnitude, specify a multipartite number

\begin{center}
$\displaystyle \overline{r_{1}r_{2} \dotsc r_{n}}$
\end{center}

\noindent associated to $m$. These multipartite numbers were first studied by MacMahon in \cite{mac}.
Let $b_{m}$ denote the number of multiplicative partitions of $m$. Note that there is a one-to-one correspondence between $b_{m}$ and the number of additive partitions of the multipartite number associated to $m$. MacMahon \cite{mac1} observed that the infinite product

$$\displaystyle \prod_{k=2}^{\infty} (1-k^{-s})^{-1}$$

\noindent is the generating function of the Dirichlet series

$$\displaystyle \sum_{m=1}^{\infty} b_{m} m^{-s}.$$

\noindent Harris and Subbarao provide a recursion for $b_{m}$ in \cite{hs} while Mattics and Dodd \cite{md} have shown that $b_{m} \leq m(\log m)^{-\alpha}$ for each fixed $\alpha >0$ and for all sufficiently large $m$. This upper bound implies a conjecture of Hughes and Shallit \cite{HS}. More precise results of an asymptotic nature on the growth rate of $b_{m}$ can be found in \cite{cep}.

In this section we consider the reciprocal Dirichlet series

$$\displaystyle \sum_{m=1}^{\infty} a_m m^{-s}$$

\noindent generated by the infinite product $\displaystyle \prod_{k=2}^{\infty} (1-k^{-s})$. The coefficients $a_{m}$ count the number of (unordered) representations of $m$ as a product of an even number of distinct integers $> 1$ minus the number of representations of $m$ as a product of an odd number of distinct integers $>1$. Note that for a positive integer $m>1$, $a_{m}$ depends only on the exponents $r_1$, $r_2$, $\dotsc$, $r_{n}$ in the canonical prime factorization of $m$. In particular, if $m$ is squarefree, the value of $a_{m}$ is a function of the number $n$ of prime factors of $m$. Let $e(n)$ denote this function. Subbarao and Verma \cite{sv} studied the asymptotic behavior of $e(n)$ and  showed that

\begin{center}
$\displaystyle \frac{\log |e(n)|}{n}$
\end{center}

\noindent is unbounded as $n \to \infty$. In fact, they prove

\begin{center}
$\displaystyle \limsup_{n \to \infty} \frac{\log |e(n)|}{n \log n} =1$.
\end{center}

\noindent Note that if we identify the factors of $m=p_{1} \dotsc p_{n}$ with subsets of $\{1, 2, \dotsc, n\}$, then $e(n)$ counts the number of ways to partition a set $S$ of $n$ elements into an even number of non-empty subsets minus the number of ways to partition $S$ into an odd number of non-empty subsets. Thus,

\begin{center}
$\displaystyle e(n)= \sum_{k=1}^{n} (-1)^{k} S(n,k)$.
\end{center}

As a result of Theorem \ref{main}, we have the following

\begin{cor} If $m$ is squarefree and contains $n$ prime factors, then $a_{m}=e(n) \neq 0$ for all $n \not \equiv 2$ and $\not\equiv 2944838 \bmod 3145728$.
\end{cor}

\subsection{$p$-adic sums}
Let $p$ be a prime. For every $a \in \mathbb{Z} \setminus \{0\}$, put

\begin{center}
$\displaystyle v_{p}(a) =$ max $\{ m \in \mathbb{Z}: p^{m} \mid a\}.$
\end{center}

\noindent We extend $v_{p}$ to $\mathbb{Q} \setminus \{0\}$ by defining
$\displaystyle v_{p}(\alpha)= v_{p}(a) - v_{p}(b)$ where $\displaystyle \alpha=\frac{a}{b}$. If we define

\[ |\alpha|_{p} =
\begin{cases}
\displaystyle p^{-v_{p}(\alpha)} \quad  \text{if} \quad \alpha \neq 0    \\
\ \\
\displaystyle 0 \quad  \text{if} \quad \alpha = 0, \\
\end{cases}
\]

\noindent then $| \quad |_{p}$ is a norm on $\mathbb{Q}$ called the {\it p-adic norm}. The field of $p$-adic numbers $\mathbb{Q}_{p}$ is the completion of $\mathbb{Q}$ with respect to $| \quad |_{p}$, i.e., $p$-adic numbers are convergent series of the form

\begin{center}
$\displaystyle \sum_{k=i}^{\infty} a_{k} p^{k},$
\end{center}

\noindent where $i$, $a_k \in \mathbb{Z}$. Recall that a $p$-adic number $\alpha \in \mathbb{Q}_{p} \setminus \mathbb{Q}$ is called a {\it p-adic irrational}.

It is a well-known result that the series $\displaystyle \sum_{n=1}^{\infty} a_{n}$ with $a_{n} \in \mathbb{Q}_{p}$ converges if and only if $|a_{n}|_{p} \to 0$ as $n \to \infty$ (see Corollary 4.1.2 in \cite{gouvea}). Thus the series

\begin{center}
$\alpha:=\displaystyle \sum_{n=1}^{\infty} n!$
\end{center}

\noindent converges in $\mathbb{Q}_{p}$ as $|n!|_{p} \to 0$. The same is true for the series

\begin{center}
$\alpha_{k}:=\displaystyle \sum_{n=1}^{\infty} n^{k} n!$
\end{center}

\noindent where $k$ is a non-negative integer. Murty and Sumner \cite{ms} investigate the irrationality of
$\alpha_{k}$. Schikhof \cite{sch} was the first to ask whether $\alpha_{0}=\alpha$ is a $p$-adic irrational or not. Murty and Sumner conjecture that it is. They also use the fact that

\begin{center}
$\displaystyle \sum_{n=0}^{m} n \cdot n! = (m+1)! - 1$
\end{center}

\noindent and $|(m+1)!|_{p} \to 0$ as $m \to \infty$ to deduce that $\alpha_{1}=-1$. Moreover, they prove using an inductive argument that

\begin{center}
$\displaystyle \alpha_{k}=v_{k} - u_{k} \alpha$,
\end{center}

\noindent where $u_{k}$, $v_{k} \in \mathbb{Z}$. In fact, they show that if one assumes that $\alpha$ is irrational, then (see Lemma 4 in \cite{ms})

\begin{center}
$\displaystyle (-1)^{k} u_{k} = \sum_{j=1}^{k+1} (-1)^{j} S(k+1, j)$.
\end{center}

\noindent As a result of this expression for $u_{k}$ and Theorem \ref{main}, we can extend Theorem 1 in \cite{ms} as follows.

\begin{cor} Let $p$ be a prime. If $\alpha$ is a $p$-adic irrational and $k+1 \not \equiv 2$ and $\not\equiv 2944838 \bmod 3145728$, then
$\alpha_{k}$ is a $p$-adic irrational.
\end{cor}

\section*{Appendix}
The following code provides the possible zeros of $f(n)$ modulo
$m$ as well as the minimal period.
\begin{verbatim}
#include <stdio.h>

#define m 'any number'

long Data1[m + 1]; long Data2[m + 1];

int main() {
        long I;
        long double Steps=0;
        Data1[1] = 1;
        Data2[1] = 0;
        for (I = 2; I < m + 1; ++I) {
                Data1[I]= 0;
                Data2[I]= 0;
        };
        long Sum=0;
        long l=0;
        printf("--------------------------------------------\n");
        printf("Possible zeros for f(n) modulo %i\n",m);
        printf("--------------------------------------------\n");
cont1:
                ++Steps;
                Sum = 0;
                for (I = m; I > 1 ; --I) {
                        Data2[I] = (m+Data1[I] * (I - 1) - Data1[I - 1]) % m;
                        Sum += Data2[I];
                };
                Data2[1] = (m - Data1[m])%m;
        if (!((Sum + Data2[1] )% m)){
                printf("Possible zero is %Lf  \n ",Steps);
        }
        // Check if minimal period is reached
        if (Sum | (Data2[1]-1)) goto Transfer;
        printf("The minimal period is : %Lf \n", Steps);
        return 0;

Transfer:
        for (I = 1; I < m + 1; ++I){
                Data1[I]=Data2[I];
        }
        goto cont1;
}
\end{verbatim}

\section*{Acknowledgments}
The authors would like to thank Ram Murty for his comments on a
preliminary version of this paper, Bruce Berndt for pointing out reference \cite{berndt}, and the referees for their encouragement and insightful comments which shortened our original proof of Theorem 1.2 and improved the exposition. The first author would also like to thank Barbara Verdonck for many productive discussions. The third author would like to mention that this paper owes its existence to a delightful talk given by Professor Murty in the Summer of 2004 at Queen's University in Kingston, Ontario, Canada.

\end{document}